\documentclass{article}
\usepackage{amsmath,amsthm,amssymb,graphicx}

\newcommand{\PP}{\mathbb{P}}
\newcommand{\ZZ}{\mathbb{Z}}
\newcommand{\QQ}{\mathbb{Q}}
\newcommand{\OPdue}{\mathcal{O}_{\PP^2}}
\newcommand{\OPtre}{\mathcal{O}_{\PP^3}}
\newcommand{\OPn}{\mathcal{O}_{\PP^n}}
\newcommand{\OPN}{\mathcal{O}_{\PP^N}}
\renewcommand{\O}{\mathcal{O}}
\newcommand{\OH}{\mathcal{O}_{H}}
\newcommand{\shI}{\mathcal{I}}
\newcommand{\shF}{\mathcal{F}}
\newcommand{\shE}{\mathcal{E}}
\newcommand{\acapo}{\hfill\break}

\newtheorem{lemma}{Lemma}[section]
\newtheorem{theorem}[lemma]{Theorem}
\newtheorem{proposition}[lemma]{Proposition}

\newtheorem*{main}{Theorem}
\newtheorem{remark}[lemma]{Remark}

\newcounter{cont}
\newenvironment{lista}
            {\begin{list}
               {\textbf{\arabic{cont}.}}
               {\setlength{\leftmargin}{16pt}
                \setlength{\rightmargin}{0pt}
                \setlength{\itemindent}{0pt}
                \setlength{\labelwidth}{10pt}
                \usecounter{cont}
               }}
              {\end{list}}

\begin{document}

\title{Non-Vanishing Theorems for Rank $2$ Bundles on $\PP^3$: a Simple Approach without the Speciality Lemma}

\author{Paolo Valabrega \and Mario Valenzano}

\date{}

\maketitle

\begin{abstract}
The paper investigates vanishing conditions on the first cohomology module of a normalized rank $2$ vector bundle $\shE$ on $\PP^3$ which force $\shE$ to split, and finds therefore strategic levels of non-vanishing for a non-split bundle. 
The present conditions improve other conditions known in the literature and are obtained with simple computations on the Euler characteristic function, avoiding the speciality lemma, Barth's restriction theorem, the discriminat property, and other heavy tools.
\end{abstract}

\section{Introduction}

The problem of finding vanishing conditions on the intermediate cohomology module that force a projective curve to be the complete intersection of two surfaces in $\PP^3$ dates back to Giuseppe Gherardelli (\cite{Ghe}), who proved the following theorem:

\begin{main}[Gherardelli]
A projective smooth curve $C$ in $\PP^3$ is the complete intersection of two surfaces if and only if the surfaces of some degree $e$ in $\PP^3$ cut out on $C$ canonical divisors (i.e.\ $C$ is $e$-subcanonical) and the linear series cut out on $C$ by the surfaces of degree $n$ is complete for all $n$ (i.e.\ $C$ is arithmetically normal).
\end{main}
 
This result is improved in \cite{CV1}, where an $e$-subcanonical (not necessarily smooth) curve C is proved to be a complete intersection under the hypothesis that $h^{1}(\shI_C(n))=0$, where $n = \frac{e}{2}+1$ if $e$ is even and $n = \frac{e+1}{2}$, or $\frac{e+3}{2}$ or $\frac{e+5}{2}$, if $e$ is odd. 
This means that there are strategic levels where the first cohomology of the ideal sheaf of a non-complete intersection curve $C$ does not vanish.

It should be observed that the problem for curves can be translated into a problem for rank $2$ vector bundles on $\PP^3$, thanks to the Serre correspondence between subcanonical curves and rank 2 bundles (see \cite{Hvb}). 
From this point of view, the Gherardelli Theorem states that a normalized non-split bundle $\shE$ has $h^1(\shE(n))\ne 0$, where $n = -1$ if $c_1 = 0$ and $n = -1$, or $0$, or $1$ if $c_1 = -1$.

By focusing on rank $2$ vector bundles, many other results were found in the years, not only in $\PP^3$ but also in $\PP^4$ or, more generally, 
$\PP^N$ or a smooth variety (see for instance \cite{EG}, \cite{KPR}, \cite{R}, \cite{RVcomp} and \cite{Valenzano}). The present paper is concerned with normalized rank 2 bundles on $\PP^3$, with Chern classes $c_2$ and $c_1 = 0$ or $-1$, and its starting points are the results contained in \cite{CV1}, \cite{CV2} and \cite{RV}. In fact \cite{RV}, improving \cite{CV1} and \cite{CV2}, shows that $h^{1}(\shE(n))=0$ is not allowed when $-r-c_1-1\leq n \leq \gamma-2$, where $\gamma$ is the level  of the \lq\lq third\rq\rq\  relevant section of the bundle and $r$ is the order of instability (see \cite{Sauer}) when the bundle is non-stable, $0$ otherwise (the lower bound actually is $-r-c_1-2$ if $\shE$ is not an instanton bundle).

The main concern of the present paper is the bound on the right side of the above inequalities (for the left side bound see \cite{CV2}, \cite{RV}). In fact our main result states that $h^1(\shE(n)) = 0$ is allowed only when $n\ge\zeta=\sqrt{3c_2+1-\frac{3c_1^2}{4}}-2-\frac{c_1}{2}$, if $\shE$ is stable, and $n > \tau = \sqrt{6\delta+1-\frac{3}{2}c_1}-2-\frac{c_1}{2}$, if $\shE$ is non-stable. So we have bounds not depending upon $\gamma$, but on the Chern classes and the degree $\delta$ of any 
minimal curve of the bundle (in the non-stable case), or (in the stable case) on  the two Chern classes  of $\shE$. We must emphasize that in the stable case the integral part of $\zeta+1$ coincides with the highest possible value for the level of the first non-zero section of the bundle (see \cite{H2} and also \cite{RVMe}). 

We observe that sometimes the present bounds are better than the known bound $\gamma-2$, while sometimes they are worse, as it can be seen in section $4$. In any event they avoid $\gamma$ and involve characters of the vector bundle that are usually easier to be computed. 

However the most remarkable feature of this paper is not the list of the results that we obtain.  What is new in this paper, in our opinion, is the very elementary approach to the problems. In fact, we want to emphasize that all our proofs avoid heavy tools like, for instance, the spectrum, the plane section and the Mumford theory of t-regular sheaves (which are used in \cite{RV}), or the reduction step (which is used in \cite{H2}). In particular they make no use of the speciality lemma, which is a key tool in \cite{CV1} (for $\PP^3$) and also in \cite{RVcomp} (for $\PP^4$). Indeed they are only based on simple techniques, involving easy numerical computations on the Euler characteristic function and the behaviour of a third degree polynomial and its roots (and, of course, on the non-elementary theory of Chern classes). We must recall that, as far as the main result of \cite{CV1} (in $\PP^3$) is concerned, \cite{P} gives a simplified proof which apparently does not involve the speciality lemma, but actually the proof is based on the following property: $c_2>0$ for a stable bundle. And it can be seen (\cite{RVcomp}, remark~8) that such a property is equivalent to the speciality lemma.

It is worth to observe that $c_2 > 0$ for a stable bundle, the so called \lq\lq discriminant property\rq\rq, is in fact a consequence of Barth's theorem on the general plane restriction of a stable bundle (see \cite{Barth}); so Barth's theorem, or at least its consequence on the \lq\lq discriminant property\rq\rq, is another heavy tool that we avoid. But of course we cannot obtain exactly the same results that require the speciality lemma, Barth's theorem and the  \lq\lq discriminant property\rq\rq. In fact our main theorems work either when the second Chern class $c_2$ is positive or when the bundle is non-stable, and these two cases cover all rank two bundles under the assumption that $c_2>0$ for a stable bundle, i.e.\ the \lq\lq discriminant property\rq\rq.

It is in particular interesting to apply our results to the main theorem of \cite{CV1}. Our simple techniques not involving the speciality lemma are sufficient to obtain \lq\lq almost\rq\rq\ such a result; in fact we obtain  all four non-vanishing theorems, but when  $c_1 = -1$, $c_2 = 2$ the non-vanishing of $h^{1}(\shE(1))$ for a non-split bundle needs non-elementary tools (see remark 3.12 and example 4.2). 

In the case of \cite{CV1} we want to emphasize that in order to exclude the case $c_2 \le 0$ and $\alpha > 0$ we need not the \lq\lq discriminant property\rq\rq\ or the speciality lemma. We are in fact able to show that the required properties hold for such a \lq\lq ghost\rq\rq\ bundle, even if we disregard the fact that Barth's theorem  and the speciality lemma exclude that it exists.

We must also observe that, while our elementary techniques work very well for the bounds $\zeta$ and $\tau$, they cannot (to our knowledge) be successfully used to find the known bound $-r-c_1-2$ (in the non-stable case).

As far as the examples of section 4 are concerned, the vanishing and non-vanishing of the cohomology have also been checked with Macaulay 2 (see \cite{Macaulay});  we wish to thank Enrico Carlini, who helped us to use it.

\section{Preliminaries} 

\begin{lista}
\item
Throughout this paper we work over an algebraically closed field $k$ of characteristic $0$. $\PP^N$ is the projective space over $k$ of dimension $N = 3$ or $2$.
\item
Let $\shE$ be a rank 2 vector bundle (i.e.\ a locally free sheaf) on $\PP^N$. We use the notation $h^i(\shE(n))$ for the dimension of the $k$-vector space $H^i(\PP^N,\shE(n))$, where $\shE(n)=\shE\otimes\OPN(n)$ for every integer $n$. 
Since $\shE^{\vee}\simeq\shE(-c_1)$, where $\shE^{\vee}$ denotes the dual of the bundle $\shE$, Serre duality on $\PP^3$ says that
$$h^i(\shE(n)) = h^{3-i}(\shE(-n-c_1-4))\qquad\text{for\ \ }i=0,\dots,3,
\text{\ \ and\ \ }\forall\, n\in\ZZ.$$
A vector bundle $\shE$ on $\PP^3$ is called ACM if it has no intermediate cohomology, i.e.\ if $h^i(\shE(n))=0$ for $i=1,2$ and for each $n\in\ZZ$.
\item
The Chern classes $c_1=c_1(\shE)$ and $c_2=c_2(\shE)$ of a rank 2 vector bundle $\shE$ will always be identified with whole numbers.
We say that $\shE$ is \emph{normalized} if $c_1$ is either $0$ or $-1$.
The Chern classes of the twisted bundle $\shE(n)$ are given by the following formulas:
\begin{align*}
c_1(\shE(n)) & = c_1 + 2n \\
c_2(\shE(n)) & = c_2 + c_1 n + n^2
\end{align*}
for all $n\in\ZZ$.
\item
Let $\shE$ be a rank 2 vector bundle on $\PP^3$ and let $H$ be a plane of $\PP^3$. Then the restriction $\shE_H=\shE\otimes\OH$ of the bundle $\shE$ to the plane $H$ is a rank 2 vector bundle on $H\cong\PP^2$, and the two bundles are linked by the following exact sequence (also called \emph{restriction sequence}):
\begin{equation}\label{restriction-sequence}
0 \to \shE(-1) \to \shE \to \shE_H \to 0.
\end{equation}
\item
Let $\shE$ be a rank 2 vector bundle on $\PP^3$ and let $n$ be any integer such that $\shE(n)$ has a non-zero section having a zero-locus $Y$ of codimension 2 (i.e.\ $Y$ is a curve). Then the vector bundle and the ideal sheaf defining $Y$ are linked by the following exact sequence (see \cite{Hvb}, theorem~1.1):
\begin{equation}
0 \to \OPtre \to \shE(n) \to \shI_Y(2n+c_1) \to 0.
\end{equation}
Moreover $\deg(Y)=c_2+c_1 n+n^2$ and $\omega_Y\simeq\O_Y(c_1-4)$. These properties hold because of the so called \emph{Serre correspondence} between rank 2 vector bundles and subcanonical curves of $\PP^3$.
\item
For every rank $2$ vector bundle $\shE$ on $\PP^3$ we define the \emph{first relevant level} of $\shE$ as
$$\alpha=\alpha(\shE):=\min\lbrace t\in\ZZ\mid h^0(\shE(t))\ne 0\rbrace,$$
the \emph{second relevant level} of $\shE$ as
$$\beta=\beta(\shE):=\min\lbrace t\in\ZZ\mid h^0(\shE(t)) > h^0(\OPn(t-\alpha))\rbrace,$$
and the \emph{third relevant level} of $\shE$ as
$$\gamma=\gamma(\shE):=\min\lbrace t\in\ZZ\mid h^0(\shE(t)) > h^0(\OPtre(t-\alpha)) + h^0(\OPtre(t-\beta))\rbrace.$$
Obviously $\alpha\le\beta\le\gamma$, and notice that $\gamma$ exists if and only if $\shE$ does not split. We denote by $\alpha_H,\beta_H,\gamma_H$ the relevant levels of $\shE_H$, with $H$ a general plane.
\item
A rank 2 vector bundle $\shE$ is a \emph{split} bundle if it is (isomorphic to) a direct sum of two line bundles, that is $\shE = \OPtre(a)\oplus\OPtre(b)$ for suitable integers $a$ and $b$. It is easy to prove that every split bundle is ACM, but also the converse is true, as stated in Horrocks' Theorem (see \cite{OSS}). This result is equivalent, in the case of rank 2 vector bundles on $\PP^3$, to the Theorem of Gherardelli, as quoted in the Introduction, through Serre correspondence cited above in no.\,5. If we disregard this result, we have in any event that: \lq\lq non-ACM\rq\rq\ implies \lq\lq non-split\rq\rq.
\item
For every rank 2 vector bundle $\shE$ on $\PP^3$ we define the number $\delta = \delta(\shE) := c_2+c_1\alpha+\alpha^2$. It is obviously true that $\delta = c_2(\shE(\alpha))$.
If the bundle $\shE$ is non-split, every non-zero section of $\shE(\alpha)$ has a zero locus of codimension 2 (see \cite{Hvb}, remark~1.0.1), and such a two-codimensional scheme associated to the bundle $\shE$ is called minimal curve of the bundle, and, by the basic properties of Chern classes, its degree is exactly the number $\delta$.
\item
A normalized rank 2 vector bundle $\shE$ on $\PP^3$ is called stable if $\alpha > 0$, semistable if $\alpha\ge-c_1$, strictly semistable if $c_1=\alpha=0$ and non-stable if $\alpha\le0$ (see \cite{H1}, lemma~3.1). The following facts are well-known when the bundle is non-stable (see \cite{Sauer}):
\begin{enumerate}
\item $\alpha_H=\alpha$,
\item $h^0(\shE(\alpha))=h^0(\shE_{H}(\alpha))=1$,
\item $\beta\ge\beta_H > -\alpha-c_1 \ge 0$,
\item $h^0(\shE(n))=h^0(\OPtre(n-\alpha)) = \displaystyle\binom{n-\alpha+3}{3}$ \ \ for\ \ $n \le -\alpha-c_1$ \ \ (really for $n\le \beta-1$),
\end{enumerate}
where $H$ is a general plane in $\PP^3$. \acapo
Also the following fact is well-known when the bundle is stable (see \cite{RV}, no.1): if $\alpha_H < \alpha$, then $\beta_H\le\alpha$.
\item
For every rank $2$ vector bundle $\shE$ on $\PP^3$ with $c_2\ge0$ we put $\zeta := \sqrt{3c_2+1-\frac{3c_1^2}{4}}-2-\frac{c_1}{2}$, that is 
$$\zeta = \sqrt{3c_2+1}-2\quad\text{if\ \ }c_1=0,\quad\text{and}\quad \zeta = \sqrt{3c_2+\frac{1}{4}}-\frac{3}{2}\quad\text{if\ \ }c_1=-1.$$
We put also $\bar\alpha := E(\zeta)+1 =$ integral part of $\zeta+1$. Observe that $\alpha\leq\bar\alpha$ by the main theorem of \cite{H2}.
\item
The Euler-Poincar\'e characteristic of a rank 2 vector bundle $\shE$ on $\PP^3$ is, by definition,
$$\chi(\shE) = h^0(\shE) - h^1(\shE) + h^2(\shE) - h^3(\shE)$$
and, thanks to the Riemann--Roch Theorem, it is a polynomial expression in the Chern classes $c_1$ and $c_2$ of $\shE$ with rational coefficients. It follows that there exists a polynomial $P(c_1,c_2;t)\in\QQ[t]$ of degree $3$, depending only on the Chern classes of $\shE$, such that
$$P(c_1,c_2;n) = \chi(\shE(n)) \qquad\forall\,n\in\ZZ.$$
Such a polynomial is also called the \emph{Hilbert polynomial} of $\shE$, and it holds that
$$P(c_1,c_2;t) = \frac{1}{3} t^3 + \left(\frac{c_1}{2} + 2\right) t^2 + 
\left(\frac{c_1^2}{2} + 2 c_1 - c_2 + \frac{11}{3}\right) t + 
\frac{c_1^3}{6} - \frac{c_1 c_2}{2} + c_1^2 + \frac{11 c_1}{6} - 2 c_2 + 2$$
i.e.
$$P(c_1,c_2;t) = \frac{1}{3}\left(t + 2 + \frac{c_1}{2}\right)\left[\left(t + 2 + \frac{c_1}{2}\right)^{\!2} - 1 + \frac{3c_1^2}{4} - 3c_2\right]\!.$$
So, for a normalized bundle, the Euler characteristic function is given by the following formulas:
\begin{align*}
& \chi(\shE(n)) = \frac{1}{3}\big(n+2\big)\Big[\big(n+2\big)^2-1-3c_2\Big] \qquad\qquad\text{if } c_1=0 \\
\noalign{\vspace{-12pt}\begin{equation}\vspace{-10pt}\label{chi}\end{equation}}
& \chi(\shE(n)) = \frac{1}{3}\left(n+\frac{3}{2}\right)\left[\left(n+\frac{3}{2}\right)^{\!2}-\frac{1}{4}-3c_2\right] \quad\,\,\,\text{if } c_1=-1.
\end{align*}

The graphical behaviour is described by the following pictures, which consider the two possible cases, i.e.\ $c_2\ge 0$ and $c_2<0$.
\vspace{.2in}
\begin{center}
\includegraphics[width=6cm]{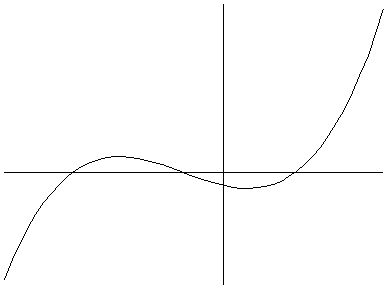}
\hspace{30pt}
\includegraphics[width=6cm]{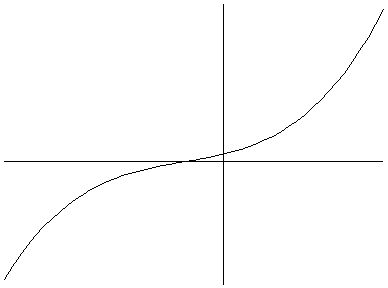}
\\
\text{Case: $c_2\ge 0$}\hspace{5cm}\text{Case: $c_2 < 0$}
\end{center}
\vspace{.2in}
In fact the Hilbert polynomial $P(c_1,c_2;t)$ of the vector bundle $\shE$, as a real function of a real variable, has three real roots if and only if $c_2 \ge 0$, while it has only one real root if and only if $c_2 < 0$.

\item
If $\shF$ is a rank 2 vector bundle on $\PP^2$, then we have:
$$\chi(\shF(n)) = (n+1)(n+2+c_1)-c_2$$
for all $n\in\ZZ$.
\item
We refer to the speciality lemma, both for varieties and for reflexive sheaves on $\PP^m$, $m\ge 3$, as it is stated and proved in \cite{RVcomp}, no.\,3 (see also \cite{Rog}). 
\acapo
Barth's theorem states that the restriction to a general plane of a stable rank 2 vector bundle on $\PP^3$  is still stable, unless it has $c_1=0$ and $c_2=1$, i.e.\ it is a null-correlation bundle (see \cite{Barth}, \cite{Hvb}, \cite{H1}).
\acapo
The \lq\lq discriminant property\rq\rq, which is a consequence of Barth's restriction theorem, states that $c_1^2 - 4 c_2 < 0$ for every stable rank 2 vector bundle on $\PP^3$. If in particular the bundle is normalized, this property states that $\alpha > 0$ implies $c_2 > 0$.
\end{lista}

\medskip
\noindent
NOTE: we agree with the following convention:
$$\binom{n}{k} =
\begin{cases}
\frac{n!}{k!\,(n-k)!} \text{\ \ if } n\ge k\\
0 \qquad\quad\text{\ \ if } n< k
\end{cases}$$

\section{Non-vanishing of the 1-cohomology and splitting of a rank 2 bundle on  $\PP^3$}

In this section we want to discuss the range of vanishing and non-vanishing of the first cohomology module of a non-split rank $2$ vector bundle on $\PP^3$. We start with the case $c_2 > 0$.
 
\begin{theorem}\label{ptre}
Let $\shE$ be a normalized non-split rank 2 vector bundle on $\PP^3$ with $c_2 > 0$. Then the following hold:
\begin{description}
\item[{}\quad\,\,$\mathrm{i)}$] $h^1(\shE(n))\ne 0$ for $-1 \le n < \zeta$.
\item[{}\quad\,$\mathrm{ii)}$] $h^1(\shE(n))\ne 0$ for $-1 \le n\le\bar\alpha-2$, and also for $n=\bar\alpha-1$ if $\zeta\notin\ZZ$.
\item[{}\quad $\mathrm{iii)}$] If $\zeta\in\ZZ$ and $\alpha < \bar\alpha$, then $h^1(\shE(\bar\alpha-1))\ne 0$.
\item[{}\quad $\mathrm{iv)}$] If $\alpha\le0$, then $h^1(\shE(\bar\alpha-1))\ne 0$.
\item[{}\quad $\mathrm{v)}$] If $\alpha>0$ and $h^1(\shE(n))=0$ with $-1\le n\le\alpha-1$, then $n=\alpha-1$ and $\alpha=\bar\alpha$.
\item[{}\quad $\mathrm{vi)}$] If $\alpha>0$ and $h^1(\shE(n))=0$ with $n\ge\alpha$, then $n\ge\bar\alpha$.
\end{description}
\end{theorem}
\begin{proof}
By the hypothesis $c_2>0$ it follows that $\zeta\ge 0$ and also $\bar\alpha-1\ge 0$.
\acapo
$\mathrm{i)}$ 
The Euler characteristic function of $\shE$ is strictly negative whenever $-1\le n<\zeta$ (see Preliminaries, no.\,11), but
$$\chi(\shE(n)) = h^0(\shE(n)) - h^1(\shE(n)) + h^2(\shE(n)) - h^3(\shE(n)),$$
where $h^2(\shE(n))\ge 0$ and $h^0(\shE(n)) - h^3(\shE(n)) = h^0(\shE(n)) - h^0(\shE(-n-4-c_1)) \ge 0$, since we have $n\ge -n-4-c_1$ for all $n\ge -1$. 
Therefore we must have $h^1(\shE(n))\ne 0$.

\smallskip\noindent
$\mathrm{ii)}$
It is a restatement of $\mathrm{i)}$ in terms of $\bar\alpha$, which is, by definition, the integral part of $\zeta+1$.

\smallskip\noindent
$\mathrm{iii)}$
If $\zeta\in\ZZ$, then $\zeta=\bar\alpha-1$, so we get
$$0 = \chi(\shE(\zeta)) = \chi(\shE(\bar\alpha-1)) = h^0(\shE(\bar\alpha-1)) - h^1(\shE(\bar\alpha-1)) + h^2(\shE(\bar\alpha-1)) - h^3(\shE(\bar\alpha-1)).$$
On the other hand $\alpha<\bar\alpha$ implies that $h^0(\shE(\bar\alpha-1))\ne 0$ and therefore $h^0(\shE(\bar\alpha-1))-h^3(\shE(\bar\alpha-1))>0$, so $\chi(\shE(\bar\alpha-1))$ cannot vanish unless $h^1(\shE(\bar\alpha-1))\ne 0$.

\smallskip\noindent
$\mathrm{iv)}$
Under the hypothesis $\alpha\le0$ we have $\alpha < \bar\alpha$, hence, if $\zeta\in\ZZ$ we apply iii), otherwise we apply ii).

\smallskip\noindent
$\mathrm{v)}$
By the hypotheses we have $\chi(\shE(n))=h^2(\shE(n))\ge 0$,; therefore $n\ge\zeta\ge\bar\alpha-1$ and it follows that $\bar\alpha-1\le n\le\alpha-1$. So $\alpha\ge\bar\alpha$ and equality holds by the main theorem of \cite{H2}. Hence  $n = \alpha-1$.

\smallskip\noindent
$\mathrm{vi)}$
In this case $\chi(\shE(n)) > 0$; therefore $n > \zeta\ge\bar\alpha-1$, so $n\ge\bar\alpha$.
\end{proof}

\begin{remark}\label{rem}
We want to emphasize that Theorem~\ref{ptre} is independent upon \cite[Theorem~7]{RVcomp}, which is the speciality lemma, and by Barth's theorem (see \cite{Barth}).  But  we cannot replace $c_2 > 0$ by $\alpha > 0$, unless we use the speciality lemma or Barth's theorem, to state the \lq\lq discriminant property\rq\rq, i.e.\ that $\alpha > 0$ implies $c_2 > 0$  (see for instance \cite[Theorem~7, Remark~8]{RVcomp} and \cite{H1}). 
\acapo
We also want to point out that, when $c_2 = 0$, $\zeta = -1$ and so \emph{i)} is an empty statement.
\end{remark}

\begin{remark}
In \cite{H2} it is proved that $\bar\alpha$ is the greatest value that $\alpha$ can reach for a rank 2 vector bundle such that $c_2\ge0$, and the proof is difficult and based on non-elementary tools like the reduction step, which necessarily involves non-locally free reflexive sheaves (outside of the category of rank two vector bundles). We want to observe that our statements, except \emph{v)}, are elementary and do not depend upon Hartshorne's theorem on the bound $\bar\alpha$. If we disregard this bound in \emph{v)} we obtain only that $\alpha\geq \bar\alpha$.
\end{remark}

\begin{remark}
In \cite{RV} it is proved (for a stable bundle) that $h^1(\shE(n)) \neq 0$ if $-1-c_1 \leq n \leq \gamma -2$ ($-2-c_1\leq n$ if $\shE$ is not an instanton bundle). The above results sometimes improve the upper inequality (see the examples below).
\end{remark}

When the bundle is non-stable we find a bound depending only upon $\delta$, as it can be seen in the theorem below. But before we state and prove it we need a few preliminary lemmas.

\begin{lemma}
The following equalities hold for every $n\in\ZZ$:
$$(n+3)(n+1) = (n+2)^2 - 1$$
and
$$(n+2)(n+1) = \left(n+\frac{3}{2}\right)^{\!2} - \frac{1}{4}.$$
As a consequence we obtain:
\begin{equation}\label{one}
\frac{1}{6}(n+3)(n+2)(n+1) = \frac{1}{6}(n+2)\Big[(n+2)^2 - 1\Big]
\end{equation}
and also
\begin{equation}\label{three}
\frac{1}{6}(n+3)(n+2)(n+1) = \frac{1}{6}\left(n+\frac{3}{2}\right)\left[\left(n+\frac{3}{2}\right)^{\!2} + 2\right] + \frac{1}{16}(4n^2+6n-1).
\end{equation}
\end{lemma}
\begin{proof}
Straightforward.
\end{proof}

\begin{lemma}\label{binomiale}
Let $\alpha$ be a fixed whole number, then for every integer $n\ge\alpha-3$ it holds
$$\binom{n-\alpha+3}{3} = \frac{1}{6}(n+3)(n+2)(n+1) - 
\left(\frac{1}{2}\alpha n^2 - \frac{1}{2}\alpha^2 n + 2\alpha n\right) - 
\left(\frac{1}{6}\alpha^3 - \alpha^2 +\frac{11}{6}\alpha\right).$$
\end{lemma}
\begin{proof}
Straightforward.
\end{proof}

\begin{lemma}\label{zerotre}
Let $\shE$ be a non-stable rank 2 vector bundle on $\PP^3$. \\
If $c_1=0$, then
$$h^0(\shE(n)) - h^3(\shE(n)) = \frac{1}{3}(n+2)\Big[(n+2)^2 - 1 + 3\alpha^2\Big]$$
for $\alpha - 3 \le n \le -\alpha - 1$, and also
$$h^0(\shE(-\alpha)) - h^3(\shE(-\alpha)) = \frac{1}{3}(-\alpha+2)\Big[(-\alpha+2)^2 - 1 + 3\alpha^2\Big] - 1.$$
If $c_1=-1$, then
$$h^0(\shE(n)) - h^3(\shE(n)) = \frac{1}{3}\left(n+\frac{3}{2}\right)\left[\left(n+\frac{3}{2}\right)^{\!2} - \frac{1}{4} + 3(\alpha^2-\alpha)\right]$$
for $\alpha - 3 \le n \le -\alpha$, and also
$$h^0(\shE(-\alpha+1)) - h^3(\shE(-\alpha+1)) = \frac{1}{3}\left(-\alpha+\frac{5}{2}\right)\left[\left(-\alpha+\frac{5}{2}\right)^{\!2} - \frac{1}{4} + 3(\alpha^2-\alpha)\right] - 1.$$
\end{lemma}
\begin{proof}
Let $c_1=0$. We have
$$h^0(\shE(n)) = \binom{n-\alpha+3}{3} \qquad\forall\,n\le-\alpha$$
and by Serre duality
$$h^3(\shE(n)) = h^0(\shE(-n-4)) = \binom{-n-\alpha-1}{3} \qquad\forall\,n\ge \alpha-4;$$
moreover it holds
$$\binom{n-\alpha+3}{3} = \frac{1}{6}(n-\alpha+3)(n-\alpha+2)(n-\alpha+1) \qquad\forall\,n\ge\alpha-3$$
and
$$\binom{-n-\alpha-1}{3} = -\frac{1}{6}(n+\alpha+1)(n+\alpha+2)(n+\alpha+3) \qquad\forall\,n\le-\alpha-1;$$
hence, for every integer $n$ such that $\alpha - 3 \le n \le -\alpha - 1$, we get by a simple computation
\begin{align*}
h^0(\shE(n)) & - h^3(\shE(n)) = \binom{n-\alpha+3}{3} - \binom{-n-\alpha-1}{3}\\
& = \frac{1}{6}\Big[(n+3-\alpha)(n+2-\alpha)(n+1-\alpha) + 
(n+3+\alpha)(n+2+\alpha)(n+1+\alpha)\Big] \\
\noalign{\vspace{4pt}}
& = \frac{1}{3}(n+3)(n+2)(n+1) + \alpha^2 (n+2) \\
\noalign{\vspace{4pt}}
& = \frac{1}{3}(n+2)\Big[(n+2)^2 - 1 + 3\alpha^2\Big].
\end{align*}
For $n=-\alpha$ we obtain on the one hand
$$h^0(\shE(-\alpha)) - h^3(\shE(-\alpha)) = h^0(\shE(-\alpha)) - 
h^0(\shE(\alpha-4)) = \binom{-2\alpha+3}{3}$$
and on the other
\[
\frac{1}{3}(-\alpha+2)\Big[(-\alpha+2)^2 - 1 + 3\alpha^2\Big] =
\frac{1}{6}(-2\alpha+3)(-2\alpha+2)(-2\alpha+1) + \frac{3\cdot2\cdot1}{6}
= \binom{-2\alpha+3}{3} + 1;
\]
therefore
$$h^0(\shE(-\alpha)) - h^3(\shE(-\alpha)) = \frac{1}{3}(-\alpha+2)\Big[(-\alpha+2)^2 - 1 + 3\alpha^2\Big] - 1.$$

\noindent
Now, let $c_1=-1$. We have
$$h^0(\shE(n)) = \binom{n-\alpha+3}{3} \qquad\forall\,n\le-\alpha+1$$
and by Serre duality
$$h^3(\shE(n)) = h^0(\shE(-n-3)) = \binom{-n-\alpha}{3} \qquad\forall\,n\ge \alpha-4;$$
moreover it holds
$$\binom{n-\alpha+3}{3} = \frac{1}{6}(n-\alpha+3)(n-\alpha+2)(n-\alpha+1) \qquad\forall\,n\ge\alpha-3$$
and
$$\binom{-n-\alpha}{3} = -\frac{1}{6}(n+\alpha)(n+\alpha+1)(n+\alpha+2) \qquad\forall\,n\le-\alpha;$$
hence, for every integer $n$ such that $\alpha - 3 \le n \le -\alpha$, we get by a simple computation
\begin{align*}
h^0(\shE(n)) & - h^3(\shE(n)) = \binom{n-\alpha+3}{3} - \binom{-n-\alpha}{3}\\
& = \frac{1}{6}\Big[(n+3-\alpha)(n+2-\alpha)(n+1-\alpha) + 
(n+2+\alpha)(n+1+\alpha)(n+\alpha)\Big] \\
\noalign{\vspace{5pt}}
& = \frac{1}{6}(2n+3)\Big[(n+2)(n+1) + 3(\alpha^2-\alpha)\Big] \\
\noalign{\vspace{5pt}}
& = \frac{1}{3}\left(n+\frac{3}{2}\right)\left[\left(n+\frac{3}{2}\right)^{\!2} - \frac{1}{4} + 3(\alpha^2-\alpha)\right]\!.
\end{align*}
For $n=-\alpha+1$ we obtain on the one hand
$$h^0(\shE(-\alpha+1)) - h^3(\shE(-\alpha+1)) = h^0(\shE(-\alpha+1)) - 
h^0(\shE(\alpha-4)) = \binom{-2\alpha+4}{3}$$
and on the other
$$\frac{1}{3}\!\!\left(-\alpha+\frac{5}{2}\right)\!\!\left[\left(-\alpha+\frac{5}{2}\right)^{\!2} - \frac{1}{4} + 3(\alpha^2-\alpha)\right]\! = 
\frac{1}{6}(-2\alpha+4)(-2\alpha+3)(-2\alpha+2) + \frac{3\cdot2\cdot1}{6} = \binom{-2\alpha+4}{3} + 1;
$$
therefore
$$h^0(\shE(-\alpha+1)) - h^3(\shE(-\alpha+1)) = \frac{1}{3}\left(-\alpha+\frac{5}{2}\right)\left[\left(-\alpha+\frac{5}{2}\right)^{\!2} - \frac{1}{4} + 3(\alpha^2-\alpha)\right] - 1.$$
\end{proof}

\begin{theorem}\label{non-stable}
Let $\shE$ be a normalized, non-split, non-stable rank 2 vector bundle on $\PP^3$. Then
\begin{description}
\item[{}\quad\,\,$\mathrm{i)}$] $h^1(\shE(n))\ne 0$ for $-1 \le n \le -\alpha-c_1$.
\item[{}\quad\,$\mathrm{ii)}$] If $\alpha=0$, then $h^1(\shE(n))\ne 0$

\bigskip
${}$\qquad for $-c_1 \le n < \tau$, where $\tau = \sqrt{6 c_2 + 1}-2$ if $c_1=0$, or

\bigskip
{}\qquad for $-c_1 \le n \le \tau$, where $\tau = \sqrt{6 c_2 + \frac{5}{2}}-\frac{3}{2}$ if $c_1=-1$.

\item[{}\quad $\mathrm{iii)}$] If $\alpha<0$, then $h^1(\shE(n))\ne 0$ for $-1 \le n \le \eta$, where 

\bigskip
${}$\qquad\qquad $\eta = \sqrt{6\delta + 1}-2$ \quad if $c_1=0$, or 

\bigskip
${}$\qquad\qquad $\eta = \sqrt{6\delta + \frac{5}{2}}-\frac{3}{2}$ \quad if $c_1=-1$.

\item[{}\quad $\mathrm{iv)}$] If moreover $c_2\geq0$ and $\alpha<0$, then $h^1(\shE(n))\ne 0$ 

\bigskip
${}$\qquad for $-\alpha-c_1 \le n < \eta$, where $\eta = \sqrt{6\delta+1-\frac{3}{4}\alpha^2}-2-\frac{3}{2}\alpha$\;\;if $c_1=0$, or

\bigskip
{}\qquad for $-\alpha-c_1 \le n \le \eta$, where $\eta = \sqrt{6\delta+\frac{13}{16}+\frac{3}{4}\alpha-\frac{3}{4}\alpha^2} - \frac{3}{4} -\frac{3}{2}\alpha$\,\;if $c_1=-1$.
\end{description}
\end{theorem}
\begin{proof}
$\mathrm{i)}$ 
Firstly we assume $c_1=0$. \acapo
Let $n$ be an integer such that $-1 \le n \le -\alpha-1$. We have, by Lemma~\ref{zerotre} and $(\ref{chi})$,
\begin{align*}
h^1(\shE(n))-h^2(\shE(n)) & = h^0(\shE(n))-h^3(\shE(n)) - \chi(\shE(n)) = \\
\noalign{\vspace{8pt}}
& = \frac{1}{3}(n+2)\Big[(n+2)^2 - 1 + 3\alpha^2\Big] - 
\frac{1}{3}(n+2)\Big[(n+2)^2 - 1 - 3 c_2\Big] = \\
\noalign{\vspace{8pt}}
& = (n+2)(c_2 + \alpha^2) > 0,
\end{align*}
since $n+2 > 0$, being $n\ge -1$, and $c_2 + \alpha^2 = \delta > 0$.
So we must have $h^1(\shE(n))\ne 0$. \\
If $n=-\alpha$, by Lemma~\ref{zerotre} we obtain
$$h^1(\shE(-\alpha))-h^2(\shE(-\alpha)) = (-\alpha+2)(c_2 + \alpha^2) - 1 > 0,$$
since $-\alpha+2\ge 2$ and $c_2 + \alpha^2 = \delta > 0$, hence $h^1(\shE(-\alpha))\ne 0$.

\medskip\noindent
Now we assume $c_1=-1$. \acapo
Let $n$ be an integer such that $-1 \le n \le -\alpha$. We have, by Lemma~\ref{zerotre} and $(\ref{chi})$,
\begin{align*}
h^1(\shE(n)) & - h^2(\shE(n)) = h^0(\shE(n))-h^3(\shE(n)) - \chi(\shE(n)) = \\
\noalign{\vspace{8pt}}
& = \frac{1}{3}\left(n+\frac{3}{2}\right)\left[\left(n+\frac{3}{2}\right)^{\!2} - \frac{1}{4} + 3(\alpha^2-\alpha)\right] - \frac{1}{3}\left(n+\frac{3}{2}\right)\left[\left(n+\frac{3}{2}\right)^{\!2} - \frac{1}{4} - 3 c_2\right] = \\
\noalign{\vspace{8pt}}
& = \left(n+\frac{3}{2}\right)(c_2-\alpha+\alpha^2) > 0,
\end{align*}
since $n+\frac{3}{2} > 0$, being $n\ge -1$, and $c_2 - \alpha + \alpha^2 = \delta > 0$.
So we must have $h^1(\shE(n))\ne 0$. \\
If $n=-\alpha+1$, by Lemma~\ref{zerotre} we obtain
$$h^1(\shE(-\alpha+1))-h^2(\shE(-\alpha+1)) = \left(-\alpha+\frac{5}{2}\right)(c_2 - \alpha + \alpha^2) - 1 > 0,$$
since $-\alpha+\frac{5}{2} > 2$ and $c_2 - \alpha + \alpha^2 = \delta > 0$, hence $h^1(\shE(-\alpha+1))\ne 0$.

\bigskip\noindent
$\mathrm{ii)}$ 
Let $c_1=0$, $n > -\alpha=0$, and assume that $h^1(\shE(n))=0$. We have
$$h^3(\shE(n)) = h^0(\shE(-n-4)) = 0 \qquad\text{and}\qquad 
h^0(\shE(n)) \ge \binom{n-\alpha+3}{3} = \binom{n+3}{3},$$
so
$$\chi(\shE(n)) - \binom{n+3}{3} \ge \chi(\shE(n)) - h^0(\shE(n)) = 
h^2(\shE(n)) \ge 0;$$
therefore by (\ref{chi}) and (\ref{one}) we can write the inequality $\chi(\shE(n)) - \binom{n+3}{3}\ge 0$ as
$$
\frac{1}{3}(n+2)\Big[(n+2)^2 - 1 - 3 c_2\Big] - \frac{1}{6}(n+2)\Big[(n+2)^2 - 1\Big] \ge 0,
$$
that is
\[
\frac{1}{6}(n+2)\Big[(n+2)^2 - 1 - 6 c_2\Big] \ge 0;
\]
but $n+2$ is strictly positive for $n>0$, so we must have
$$n \ge \sqrt{6 c_2 + 1} -2.$$

\medskip\noindent
Now, let $c_1=-1$, $n > -\alpha+1=1$, and assume that $h^1(\shE(n))=0$. We have
$$h^3(\shE(n)) = h^0(\shE(-n-3)) = 0 \qquad\text{and}\qquad 
h^0(\shE(n)) \ge \binom{n-\alpha+3}{3} = \binom{n+3}{3},$$
so
$$\chi(\shE(n)) - \binom{n+3}{3} \ge \chi(\shE(n)) - h^0(\shE(n)) = 
h^2(\shE(n)) \ge 0.$$
By (\ref{chi}) and (\ref{three}) we can write the inequality $\chi(\shE(n)) - \binom{n+3}{3}\ge 0$ as
\[
\frac{1}{3}\left(n+\frac{3}{2}\right)\left[\left(n+\frac{3}{2}\right)^{\!2} - \frac{1}{4} - 3 c_2\right] - \frac{1}{6}\left(n+\frac{3}{2}\right)\left[\left(n+\frac{3}{2}\right)^{\!2} + 2\right] - \frac{1}{16}(4n^2+6n-1) \ge 0,
\]
that is
\[
\frac{1}{6}\left(n+\frac{3}{2}\right)\left[\left(n+\frac{3}{2}\right)^{\!2} - \frac{5}{2} - 6 c_2\right] - \frac{1}{16}(4n^2+6n-1) \ge 0;
\]
notice now that
$$-\frac{1}{16}(4n^2+6n-1) < 0 \qquad\forall\,n>1;$$
therefore we obtain
\[
\frac{1}{6}\left(n+\frac{3}{2}\right)\left[\left(n+\frac{3}{2}\right)^{\!2} - \frac{5}{2} - 6 c_2\right] > 0,
\]
but $n+\frac{3}{2}$ is strictly positive for $n>1$, so we must have
$$n > \sqrt{6 c_2 + \frac{5}{2}} - \frac{3}{2}.$$

\smallskip\noindent
$\mathrm{iii)}$
By i) we can assume that $\eta\ge-\alpha-c_1$. \acapo
Firstly we consider the case $c_1=0$ and $\alpha<0$.
Let $n > -\alpha$ and assume that $h^1(\shE(n))=0$. We have
$$h^3(\shE(n)) = h^0(\shE(-n-4)) = 0 \qquad\text{and}\qquad 
h^0(\shE(n)) \ge \binom{n-\alpha+3}{3},$$
so
$$\chi(\shE(n)) - \binom{n-\alpha+3}{3} \ge \chi(\shE(n)) - h^0(\shE(n)) = 
h^2(\shE(n)) \ge 0;$$
therefore by (\ref{chi}), Lemma~\ref{binomiale} and (\ref{one}) we can write the inequality $\chi(\shE(n)) - \binom{n-\alpha+3}{3}\ge 0$ as
$$
\frac{1}{3}(n+2)\Big[(n+2)^2 - 1 - 3 c_2\Big] -\frac{1}{6}(n+2)\Big[(n+2)^2 - 1\Big] + \left(\frac{1}{2}\alpha n^2 - \frac{1}{2}\alpha^2 n + 2\alpha n\right) + \left(\frac{1}{6}\alpha^3 - \alpha^2 +\frac{11}{6}\alpha\right) \ge 0,
$$
hence, by adding and subtracting the quantity $\alpha^2(n+2)$, we get
\begin{align*}
\frac{1}{6}(n+2)\Big[(n+2)^2 - 1 - 6 c_2\Big] & + \frac{1}{6}(n+2)(-6\alpha^2)
+ \alpha^2 n + 2\alpha^2 + {} \\
\noalign{\vspace{4pt}}
& + \left(\frac{1}{2}\alpha n^2 - \frac{1}{2}\alpha^2 n + 2\alpha n\right) + \left(\frac{1}{6}\alpha^3 - \alpha^2 +\frac{11}{6}\alpha\right) \ge 0,
\end{align*}
that is
$$
\frac{1}{6}(n+2)\Big[(n+2)^2 - 1 - 6 (c_2 + \alpha^2)\Big] + \left(\frac{1}{2}\alpha n^2 + \frac{1}{2}\alpha^2 n + 2\alpha n\right) + 2\alpha^2 + 
\left(\frac{1}{6}\alpha^3 - \alpha^2 +\frac{11}{6}\alpha\right) \ge 0,
$$
which we can write as
$$
\frac{1}{6}(n+2)\Big[(n+2)^2 - 1 - 6 \delta\Big] + \frac{1}{2}\alpha n(n+\alpha+4) + 2\alpha^2 + \frac{1}{6}\alpha(\alpha^2 - 6\alpha + 11) \ge 0;
$$
notice now that
$$n > -\alpha,\quad n+\alpha+4 > 4 \text{\quad and\quad} \alpha<0,$$
so it holds
$$\frac{1}{2}\alpha n(n+\alpha+4) < -2\alpha^2,$$
i.e.\
$$\frac{1}{2}\alpha n(n+\alpha+4) + 2\alpha^2 < 0$$
for all $n > -\alpha$; moreover
$$\frac{1}{6}\alpha(\alpha^2 - 6\alpha + 11) < 0 \text{\qquad for all } \alpha < 0$$
since $\alpha^2 - 6\alpha + 11 > 0$ for all $\alpha$ (the discriminant of this quadratic form is strictly negative); therefore we obtain
$$\frac{1}{6}(n+2)\Big[(n+2)^2 - 1 - 6 \delta\Big] > 0,$$
but $n+2$ is strictly positive for $n>-\alpha\ge 1$, so we must have
$$n > \sqrt{6\delta + 1} -2.$$

\medskip\noindent
Now we consider the case $c_1=-1$ and $\alpha<0$.
Let $n > -\alpha+1$ and assume that $h^1(\shE(n))=0$. We have
$$h^3(\shE(n)) = h^0(\shE(-n-3)) = 0 \qquad\text{and}\qquad 
h^0(\shE(n)) \ge \binom{n-\alpha+3}{3},$$
so
$$\chi(\shE(n)) - \binom{n-\alpha+3}{3} \ge \chi(\shE(n)) - h^0(\shE(n)) = 
h^2(\shE(n)) \ge 0.$$
By (\ref{chi}), Lemma~\ref{binomiale} and (\ref{three}) we can write the inequality $\chi(\shE(n)) - \binom{n-\alpha+3}{3} \ge 0$ as
\begin{align*}
\frac{1}{3}\left(n+\frac{3}{2}\right)\left[\left(n+\frac{3}{2}\right)^{\!2} - \frac{1}{4} - 3 c_2\right]  & - \frac{1}{6}\left(n+\frac{3}{2}\right)\left[\left(n+\frac{3}{2}\right)^{\!2} + 2\right] - \frac{1}{16}(4n^2+6n-1) + {} \\
\noalign{\vspace{4pt}}
& + \left(\frac{1}{2}\alpha n^2 - \frac{1}{2}\alpha^2 n + 2\alpha n\right) + \left(\frac{1}{6}\alpha^3 - \alpha^2 +\frac{11}{6}\alpha\right) \ge 0
\end{align*}
hence, by adding and subtracting the quantity $(\alpha^2-\alpha)\left(n+\frac{3}{2}\right)$, we get
\begin{align*}
\frac{1}{6}\left(n+\frac{3}{2}\right)\Bigg[\!\left(n+\frac{3}{2}\right)^{\!2} & - \frac{1}{4} - 6 c_2\Bigg] - \frac{1}{6}\left(n+\frac{3}{2}\right)\left[\left(n+\frac{3}{2}\right)^{\!2} + 2\right] - \frac{1}{16}(4n^2+6n-1) + {} \\
\noalign{\vspace{4pt}}
& + \frac{1}{6}\left(n+\frac{3}{2}\right)\big(\!-6(\alpha^2-\alpha)\big) + (\alpha^2 n -\alpha n) + \frac{3}{2}(\alpha^2 -\alpha) + {} \\
\noalign{\vspace{4pt}}
& + \left(\frac{1}{2}\alpha n^2 - \frac{1}{2}\alpha^2 n + 2\alpha n\right) + \left(\frac{1}{6}\alpha^3 - \alpha^2 +\frac{11}{6}\alpha\right) \ge 0,
\end{align*}
that is
\begin{align*}
\frac{1}{6}\left(n+\frac{3}{2}\right) &\left[\left(n+\frac{3}{2}\right)^{\!2} - \frac{1}{2} - 2 - 6(c_2-\alpha+\alpha^2)\right] - \frac{1}{16}(4n^2+6n-1) +{} \\
\noalign{\vspace{4pt}}
& + \left(\frac{1}{2}\alpha n^2 + \frac{1}{2}\alpha^2 n + \alpha n\right) + \frac{3}{2}(\alpha^2 - \alpha) +
\left(\frac{1}{6}\alpha^3 - \alpha^2 +\frac{11}{6}\alpha\right) \ge 0,
\end{align*}
which we can write as
\begin{align*}
\frac{1}{6}\left(n+\frac{3}{2}\right) &\left[\left(n+\frac{3}{2}\right)^{\!2} - \frac{5}{2} - 6\delta\right] - \frac{1}{16}(4n^2+6n-1) + {} \\
\noalign{\vspace{4pt}}
& + \frac{1}{2}\alpha n(n+\alpha+2) + \frac{3}{2}(\alpha^2 - \alpha) + \frac{1}{6}\alpha(\alpha^2 - 6\alpha + 11) \ge 0;
\end{align*}
notice now that
\[
-\frac{1}{16}(4n^2+6n-1) < 0 \qquad \forall\,n > -\alpha+1 \ge 2
\]
(in fact we have $-\frac{1}{16}(4n^2+6n-1) < 0$ for all $n \ge 1$), and also
\[
\frac{1}{2}\alpha n(n+\alpha+2) + \frac{3}{2}(\alpha^2 - \alpha) < 0 \qquad \forall\,n > -\alpha+1 
\]
\[
\frac{1}{6}\alpha(\alpha^2 - 6\alpha + 11) < 0 \qquad\forall\,\alpha < 0;
\]
therefore we obtain
$$\frac{1}{6}\left(n+\frac{3}{2}\right)\left[\left(n+\frac{3}{2}\right)^{\!2} - \frac{5}{2} - 6\delta\right] > 0,$$
but $n+\frac{3}{2}$ is strictly positive for $n>-\alpha+1\ge 2$, so we must have
$$n > \sqrt{6\delta + \frac{5}{2}} - \frac{3}{2}.$$

\smallskip\noindent
$\mathrm{iv)}$ 
Let $n > -\alpha-c_1$ and assume that $h^1(\shE(n))=0$. We have
$$h^3(\shE(n)) = h^0(\shE(-n-4-c_1)) = 0 \qquad\text{and}\qquad 
h^0(\shE(n)) \ge \binom{n-\alpha+3}{3},$$
so
$$\chi(\shE(n)) - \binom{n-\alpha+3}{3} \ge \chi(\shE(n)) - h^0(\shE(n)) = 
h^2(\shE(n)) \ge 0.$$
Firstly we consider the case $c_1=0$. 
Since the quantity $\frac{1}{6}(-\alpha^3+\alpha)$ is positive or 0, we obtain the following inequality
$$\frac{1}{3}(n+2)\Big[(n+2)^2 - 1 - 3 c_2\Big] - \binom{n-\alpha+3}{3} + \frac{1}{6}(-\alpha^3+\alpha) \ge 0,$$
which we can write as
$$\frac{1}{6}(n+2)\Bigg[\!\left(n+2+\frac{3}{2}\alpha\right)^{\!2} + \frac{3}{4}\alpha^2 - 1 - 6 \delta\Bigg] \ge 0;$$
but $6\delta+1-\frac{3}{4}\alpha^2=6 c_2+\frac{21}{4}\alpha^2+1$ is positive, because $c_2\geq0$; therefore we must have
$$n\ge \sqrt{6\delta+1-\frac{3}{4}\alpha^2}-2-\frac{3}{2}\alpha.$$
Now we consider the case $c_1=-1$. 
Since the quantity $\frac{1}{6}\!\left(-\alpha+\frac{3}{2}\right)\left(\alpha^2-\frac{1}{4}\right)$ is strictly positive, we obtain the following inequality
$$\frac{1}{3}\left(n+\frac{3}{2}\right)\Bigg[\!\left(n+\frac{3}{2}\right)^2 - 1 - 3 c_2\Bigg] - \binom{n-\alpha+3}{3} + \frac{1}{6}\left(-\alpha+\frac{3}{2}\right)\left(\alpha^2-\frac{1}{4}\right) > 0,$$
which we can write as
$$\frac{1}{6}\left(n+\frac{3}{2}\right)\Bigg[\!\left(n+\frac{3}{4}+\frac{3}{2}\alpha\right)^{\!2} - \frac{13}{16} - \frac{3}{4}\alpha + \frac{3}{4}\alpha^2 - 6 \delta\Bigg] > 0;$$
but $6\delta+\frac{13}{16}+\frac{3}{4}\alpha-\frac{3}{4}\alpha^2=6 c_2+\frac{21}{4}(\alpha^2-\alpha)+\frac{13}{16}$ is strictly positive, because $c_2\geq0$; therefore we must have
$$n > \sqrt{6\delta+\frac{13}{16}+\frac{3}{4}\alpha-\frac{3}{4}\alpha^2} - \frac{3}{4} -\frac{3}{2}\alpha.$$
\end{proof}

\begin{remark}
The above theorem includes as a special case the bundles which are strictly semistable $(c_1 = \alpha = 0)$.
\end{remark}

\begin{remark}
In \cite{RV} it is proved that $h^1(\shE(n)) \neq 0$ if $-r-c_1-2 \leq n \leq \gamma -2$ (where $r$ is the order of instability of $\shE$, as defined in \cite{Sauer}, i.e.\ $r = -\alpha-c_1$ if the bundle is non-stable and $r=0$ if the bundle is stable). The above result sometimes improves the upper inequality (see the examples below).
\end{remark}

\begin{remark}
In order to obtain the lower bound $-r-c_1-2=\alpha-2$ in the non-stable case, our direct technique, based on the study of the Hilbert polynomial of the bundle, seems to be of no help and, to our knowledge, it is necessary to make use of less elementary tools (see \cite{CV2} and \cite{RV}).
\end{remark}

\begin{remark}
As a consequence of the above theorems we \lq\lq almost\rq\rq\ obtain the main theorem of \cite{CV1}, with the only exception $h^1(\shE(1))\ne0$ for a stable non-split bundle $\shE$ with $c_1=-1$ and $c_2=2$.
This case requires (to our knowledge) the use of Barth's and Castelnuovo-Mumford's theorems (see \cite{HS} and example 4.2 in section 4). 
\acapo 
We want to observe that the case $\alpha > 0$, $c_2 \leq 0$ is excluded by Barth's theorem (see remark~\ref{rem} above), but our proof of the main result of \cite{CV1} does not need such a strong theorem, we see directly that in the (impossible) event of a stable vector bundle with negative $c_2$ the main claim of \cite{CV1} must hold (we really give two proofs, both independent upon Barth's theorem -- see subcase~2.5 below).
\end{remark}

For our proof we need two lemmas. 
 
\begin{lemma}\label{split}
Let $\shE$ be a rank 2 vector bundle on $\PP^N$, $N\ge 2$. Then $\shE$ splits if and only if $\delta=0$.
\end{lemma}
\begin{proof}
If $\shE$ splits, than $\shE = \O(a)\oplus\O(-a+c_1)$, for some non-negative integer $a$. Therefore $\delta$ is easily seen to be $0$.
Assume now that $\shE$ is a non-split bundle having $\delta = 0$. Then (see  Preliminaries, no.\,8), $\shE(\alpha)$ has a non-zero section giving rise to a two-codimensional scheme, whose degree, by \cite{AG}, Appendix A, 3, C6, is exactly $\delta$. Hence we obtain a contradiction.
\end{proof}

\begin{lemma}\label{left-vanishing}
Let $\shE$ be a rank 2 vector bundle on $\PP^3$ with first relevant level $\alpha$. If $h^1(\shE(m))=0$ with $m\le\alpha-2$, then $h^1(\shE(n))=0$ for all $n\le m$.
\end{lemma}
\begin{proof}
Let $H$ be a general plane and consider the restriction sequence $(\ref{restriction-sequence})$ twisted by $m+1$. In cohomology we obtain
$$0=H^0(\shE(m+1)) \to H^0(\shE_H(m+1)) \to H^1(\shE(m))=0$$
since $m+1\le\alpha-1$. So we get $h^0(\shE_H(n))=0$ for all $n\le m+1$, which implies $h^1(\shE(n-1)) \le h^1(\shE(n))$ for all $n\le m$; but $h^1(\shE(m))=0$, hence the thesis.
\end{proof}

We are now able to state and prove the following splitting theorem, which is \lq\lq almost\rq\rq\ the main theorem of \cite{CV1} (where this splitting criterion is stated and proved for subcanonical curves of $\PP^3$).

\begin{proposition}
Let $\shE$ be a normalized rank 2 vector bundle on $\PP^3$. Then $\shE$ is a split bundle if and only if
\begin{description}
\item[{}\quad\,\,$\mathrm{a)}$] $h^1(\shE(-1))=0$ if $c_1=0$,
\item[{}\quad\,\,$\mathrm{b)}$] $h^1(\shE(-1))=0$ or $h^1(\shE)=0$ if $c_1=-1$.

If moreover $c_1 = -1$ but $\shE$ is not a stable bundle with $c_2 = 2$, then it splits if and only if
\item[{}\quad\,\,$\mathrm{c)}$] $h^1(\shE(1))=0$ if $c_1=-1$.
\end{description}
\end{proposition}
\begin{proof}
If the bundle is split, then it is ACM. So we assume that $\shE$ is non-split. Recall that if $c_1 = -1$ then $c_2$ must be even. We distinguish two cases (stable and non-stable) and a few subcases.
\acapo
\underline{Case 1}: $\shE$ is non-stable, i.e.\ $\alpha \leq 0$. \acapo
By Theorem~\ref{non-stable}, i), we know that $ h^1(\shE(n))\ne 0$ for every $n$ such that $-1 \le n \le -\alpha-c_1$, where $-\alpha-c_1 \ge 0$ if $c_1 = 0$,  $-\alpha-c_1 \ge 1$ if $c_1 = -1$.  Therefore our four claims follow immediately.
\acapo
\underline{Case 2}: $\shE$ is stable, i.e.\ $\alpha > 0$.
\acapo
\underline{Subcase 2.1}: $c_1=0$ and $c_2 > 0$. 
Since $\zeta\ge 0$, by Theorem~\ref{ptre}, i), it results for sure $h^1(\shE(-1))\ne 0$.
\acapo
\underline{Subcase 2.2}: $c_1=-1$ and $c_2\ge 4$. 
Since $\zeta\ge 2$, by Theorem~\ref{ptre}, i), it holds $h^1(\shE(n))\ne 0$ for $n=-1,0,1$. 
\acapo
\underline{Subcase 2.3}: $c_1=-1$ and $c_2 = 2$. 
In this event we have $\zeta = 1$ and so Theorem~\ref{ptre}, i), implies that $ h^1(\shE(-1))\ne 0$ and $ h^1(\shE)\ne 0$, but not $ h^1(\shE(1))\ne 0$.
\acapo
\underline{Subcase 2.4}: $c_1=0$ and $c_2\le 0$ (proof independent upon Barth's theorem). Since $\alpha > 0$ and $c_2\le 0$, it holds that $-h^1(\shE)+h^2(\shE)=\chi(\shE)>0$, so we must have $h^1(\shE(-4))=h^2(\shE)\ne 0$. Since $\alpha-2\ge-1$, it follows by Lemma~\ref{left-vanishing} that $h^1(\shE(-1))\ne 0$.
\acapo
\underline{Subcase 2.5}: $c_1=-1$ and $c_2\le 0$. \acapo
(First proof independent upon Barth's theorem). \acapo
If $ h^1(\shE(-1)) = 0$ by Lemma~\ref{left-vanishing} we also have $h^1(\shE(-2)) = 0$ and moreover $h^0(\shE(-1)) = 0$ by hypothesis and so also $h^3(\shE(-1)) = 0$. Therefore $0 = \chi(\shE(-1)) = -\frac{c_2}{2}$. Hence we have: $c_2 = 0$. But $- h^1(\shE)+ h^2(\shE) = \chi(\shE) = 1-\frac{3c_2}{2} > 0$, and this implies $h^2(\shE) =  h^1(\shE(-3)) \ne 0$, absurd by Lemma~\ref{left-vanishing}. Therefore we can conclude that $ h^1(\shE(-1))\ne 0$ in any event.
\acapo
Let us now assume that $h^1(\shE) = 0$. 
First of all observe that $\delta = c_2-\alpha+\alpha^2 \ge 1$ (Lemma~\ref{split}). Since by Lemma~\ref{left-vanishing} we can assume that $\alpha = 1$, we have: $\delta = c_2 \ge 1$. Hence  $c_2$  must be positive if $h^1(\shE) = 0$, and this contradicts the hypothesis. 
\acapo
Let us now assume that $h^1(\shE(1)) = 0$.
By Lemma~\ref{left-vanishing} $\alpha$ must be at most $2$. If $\alpha = 1$ we obtain that $\delta = c_2 \ge 1$ (see Lemma~\ref{split}), so $\alpha = 2$ exactly. Since $\delta = c_2-2+4 = c_2+2$, $c_2$ cannot be strictly negative, so it is $0$. 
Let us now consider the number $\alpha_H$, level of the first relevant section of the general plane restriction $\shE_H$. Then we have the following exact sequences:
\begin{equation*}
0 \to \OPtre \to \shE(2) \to \shI_Y(3) \to 0
\end{equation*}
\begin{equation*}
0 \to \OPdue \to \shE_H(2) \to \shI_{Y\cap H}(3) \to 0.
\end{equation*}
Therefore we obtain that $\alpha_H$ cannot be less than $0$, so it is $0$, $1$ or $2$. Assume that it is $2$. Then we have $4=\chi(\shE_H(1)) = - h^1\shE_H(1)$, which is absurd. Hence $\alpha_H$ is either $0$ or $1$. But in both events we obtain a contradiction, because $\delta_H = c_2-\alpha_H+\alpha_H^2$ cannot be $0$ (Lemma~\ref{split}). 
We conclude that a vector bundle $\shE$ with $c_1 = -1$, $c_2 \leq 0$, $\alpha > 0$,  $h^1(\shE(1)) = 0$ cannot exist. \acapo
(Second proof independent upon Barth's theorem). \acapo
Since $\alpha > 0$ and $c_2\le 0$, it holds that $-h^1(\shE)+h^2(\shE)=\chi(\shE)>0$, so we must have $h^1(\shE(-3))=h^2(\shE)\ne 0$. Since $\alpha-2\ge-1$, it follows by Lemma~\ref{left-vanishing} that $h^1(\shE(n))\ne 0$ for $-3\le n\le\alpha-2$.
If $\alpha=1$, then $c_2 = c_2-\alpha+\alpha^2 = \delta > 0$, since $\shE$ is non-split, but this contradicts the assumption $c_2\le 0$.
Therefore we must have $\alpha\ge 2$, which implies $h^1(\shE(-1))\ne 0$ and $h^1(\shE)\ne 0$. Moreover, if $\alpha\ge 3$ we also have $h^1(\shE(1))\ne 0$.
If $\alpha=2$, then $c_2 + 2 = c_2 -\alpha+\alpha^2=\delta>0$, since $\shE$ is non-split, so we must have $c_2>-2$, which implies $c_2=0$, since $c_2\le 0$ and moreover $c_2$ is even. So, let us now assume that $c_1=-1$, $c_2=0$, $\alpha=2$ and $h^1(\shE(1)) = 0$. Let $H$ be a general plane, then we have $\chi(\shE_H(n)) = (n+1)^2 >0$ for all $n\ge 0$. Therefore $h^0(\shE_H)\ne 0$, i.e.\ $\alpha_H\le 0$. Moreover $\alpha_H^2-\alpha_H=\delta(\shE_H)>0$ by the non-splitting hypothesis, so $\alpha_H\le-1$, which implies that $\beta_H>-\alpha_H-c_1\ge 2=\alpha$, in contradiction with $\beta_H\le\alpha$, which is a consequence of $\alpha_H < \alpha$ (see Preliminaries~no.\,9). Therefore we must have $h^1(\shE(1))\ne 0$ also when $\alpha=2$.
\end{proof}

\section{Examples of rank 2 vector bundles on $\PP^3$}

\subsection{Stable bundles with $c_1=0$, $c_2=2$.}
Such a bundle has $\alpha =1$ and, moreover, $h^1(\shE)\neq 0,h^1(\shE(1)) = 0$ (see \cite{Hvb}, no.\,9). This agrees with Theorem \ref{ptre}.

\subsection{Stable bundles with $c_1=-1$, $c_2=2$.}
The minimal curve of such a bundle is the disjoint union of two irreducible conics and $\zeta=1\in\ZZ$, $\bar\alpha=2$ (see \cite{HS}); the cohomology of the bundle is described by the following table:
\begin{center}
\renewcommand{\arraystretch}{1.25}
\begin{tabular}{|c|r r r r r r|}
\hline
$h^1$ & 0 & 1 & 2 & 1 & 0 & 0 \\
$h^0$ & 0 & 0 & 0 & 1 & 7 & 21 \\
\hline
& $-2$ & $-1$ & 0 & 1 & 2 & 3 \\
\hline
\end{tabular}
\end{center}
In this case $\bar\alpha$  equals $\gamma$; since $\zeta$ is an integer, we see that $h^1(\shE)\ne 0$, by Theorem~\ref{ptre}, ii), which gives the same result as \cite{RV}.  But, since $\alpha = 1 < \bar\alpha = 2$, we can also apply Theorem~\ref{ptre}, iii), and see that $h^1(\shE(1))\ne 0$. Therefore Theorem~\ref{ptre} gives exactly the highest non-vanishing 1-cohomology.
\acapo
Observe that $h^1\shE(1)) \neq 0$ when $c_1 = -1, c_2 = 2$ is exactly the missing part of the splitting criterion of \cite{CV1}.

\subsection{Stable bundles with $c_1=0$, $c_2=4$.}
We have $\zeta=\sqrt{13}-2\notin\ZZ$ and $\bar\alpha=2$, and we must distinguish 3 cases (see \cite{Chang}). \acapo
Case A: the minimal curve of $\shE$ is the disjoint union of two elliptic quartics, and the cohomology of the bundle is described by the following table:
\begin{center}
\renewcommand{\arraystretch}{1.25}
\begin{tabular}{|c|r r r r r r r|}
\hline
$h^1$ & 0 & 1 & 4 & 6 & 4 & 1 & 0 \\
$h^0$ & 0 & 0 & 0 & 0 & 0 & 5 & 20 \\
\hline
& $-3$ & $-2$ & $-1$ & 0 & 1 & 2 & 3 \\
\hline
\end{tabular}
\end{center}
Observe that in this case $\bar\alpha-1=1$ and so $h^1(\shE(1))\neq 0$. Therefore Theorem~\ref{ptre} gives a better description of the cohomology than \cite{RV}, because $\gamma-2=0$.
\acapo
Case B: the minimal curve of the bundle $\shE$ is the disjoint union of an elliptic cubic and an elliptic quintic, and the cohomology of the bundle is described by the following table:
\begin{center}
\renewcommand{\arraystretch}{1.25}
\begin{tabular}{|c|r r r r r r r|}
\hline
$h^1$ & 0 & 1 & 4 & 6 & 4 & 2 & 0 \\
$h^0$ & 0 & 0 & 0 & 0 & 0 & 6 & 20 \\
\hline
& $-3$ & $-2$ & $-1$ & 0 & 1 & 2 & 3 \\
\hline
\end{tabular}
\end{center}
Observe that $\bar\alpha-1=1> \gamma-2=0$ and Theorem~\ref{ptre} is better than \cite{RV}. \acapo
Case C: the minimal curve of the bundle $\shE$ is the disjoint union of a line and a double conic, and the cohomology of the bundle is described by the following table:
\begin{center}
\renewcommand{\arraystretch}{1.25}
\begin{tabular}{|c|r r r r r r r|}
\hline
$h^1$ & 0 & 1 & 4 & 6 & 5 & 2 & 0 \\
$h^0$ & 0 & 0 & 0 & 0 & 1 & 6 & 20 \\
\hline
& $-3$ & $-2$ & $-1$ & 0 & 1 & 2 & 3 \\
\hline
\end{tabular}
\end{center}
Observe again that $\bar\alpha-1=1$ while $\gamma-2=0$.
\acapo
Observe that, in all cases, $h^1(\shE(-3)) = 0$, hence the examples are sharp for the lower bound of \cite{CV2} and \cite{RV} (semistable case with $r = 0$).

\subsection{Bundles with natural cohomology.}
It is well-known (see \cite{Hi}) that, with few exceptions, there are rank two stable vector bundles with given Chern classes and natural cohomology. For such a bundle $\alpha = \bar\alpha$ (see \cite{RVMe}, esempi, (i)). Now choose $c_1 = 0$ and $c_2$ in such a way that $\zeta$ is not an integer and, moreover,  $\alpha\geq7$. Since we have $h^0(\shE(\alpha)) = \frac{1}{3}(\alpha+2)[(\alpha+2)^2-1-3c_2] \geq1$, it is clear that $\alpha \geq7$ implies $h^0(\shE(\alpha))\geq3$, hence $\bar\alpha = \alpha =\beta = \gamma$ and so $\bar\alpha -1 > \gamma - 2$. 
\acapo
By Theorem~\ref{ptre}, ii), $h^1(\shE(\bar\alpha-1))\neq 0$, while \cite{RV} gives $h^1(\shE(\gamma-2))\neq 0$.

\subsection{A non-stable bundle with $c_1 = 0$, $c_2=9$, $\alpha = -3$.}
Take in $\PP^3$ homogeneous coordinates $x,y,z,t$ and let $Y$ be the non-reduced structure on the line $L\colon x=y=0$ defined by the ideal $I=(x^{18},x^{15}y,x^{12}y^2,x^9y^3,x^6y^4,x^3y^5,y^6,z^3x^3-yt^5)$ (see \cite{CV1}, example~3.1, (iii)). Then $Y$ is $(-10)$-subcanonical and the zero locus of a section of a bundle $\shF = \shE(-3)$, where $\shE$ has $c_1 = 0$, $\alpha = -3$, $\gamma = 9$, $c_2 = 9$, $\delta = 18$.  By Theorem~\ref{non-stable}, iii), we see that $h^1(\shE(8))\ne0$, while the fact that $\gamma-2 = 7$ implies only that $h^1(\shE(7))\ne0$. But if we apply Theorem~\ref{non-stable}, iv), which is possible because $c_2\ge0$, we obtain a much better result: $h^1(\shE(12))\ne0$.

\subsection{A strictly semistable bundle with $c_1 = \alpha = 0$, $c_2=3$.}
Take in $\PP^3$ homogeneous coordinates $x,y,z,t$ and let $Y$ be the non-reduced structure on the line $L\colon x=y=0$ defined by the homogeneous ideal $I = (x^3,x^2y,xy^2,y^3,z^2x-yt^2)$ (see \cite{CV1}, example~3.1). Then $Y$ is $(-4)$-subcanonical and the zero locus of a section of a bundle $\shE$ having $c_1 = 0$, $\alpha = 0$, $\beta = \gamma = 3$, $\delta = c_2 = 3$. Theorem~\ref{non-stable}, ii), implies that $h^1(\shE(2))\ne0$, while $h^1(\shE(1))\ne0$ is the best we can deduce from the value $1$ of $\gamma-2$. We must observe that $h^1(\shE(3)) = 0$ (see \cite{CV1}, example~3.1 and \cite{C}, example~3.4). 
In this case also Theorem~\ref{ptre} can be applied, but it gives a worse bound ($h^1(\shE(1))\ne0$).

\subsection{A stable bundle with $c_1 = 0$, $c_2=47$, $\alpha = 1$.}
Take in $\PP^3$ homogeneous coordinates $x,y,z,t$ and let $Y$ be the non-reduced structure on the line $L\colon x=y=0$ defined by the homogeneous ideal $I = (x^{16},x^{12}y^8,x^8y^{16},x^4y^{24},y^{32},z^6x^4-y^8t^2)$ (see \cite{CV1}, example~3.1). Then $Y$ is $14$-subcanonical and the zero locus of a section of a bundle $\shF = \shE(9)$, where $\shE$ has  $c_1 = 0$, $\alpha = 1$, $\delta = 48$, $c_2 = 47$. Then Theorem~\ref{ptre}, ii), gives $h^1(\shE(9))\ne0$. We must observe that in this case $\gamma -2 = 7$.

\subsection{A stable bundle with $c_1 = 0$, $c_2=20$, $\alpha = 2$, $\gamma = 10$.}
Take in $\PP^3$ homogeneous coordinates $x,y,z,t$ and let $Y$ be the non-reduced structure on the line $L\colon x=y=0$ defined by the homogeneous ideal $I = (x^{12},x^{10}y^2,x^8y^{4},x^6y^{6},x^{4}y^{8},x^{2}y^{10},$ $y^{12},z^2x^2-y^2t^2)$ (see \cite{CV1}, example~3.1). Then $Y$ is $0$-subcanonical and the zero locus of a section of a bundle $\shF = \shE(2)$, where $\shE$ has  $c_1 = 0$, $\alpha = 2$, $\gamma = 10$, $c_2 = 20$. Then Theorem~\ref{ptre}, ii), gives $h^1(\shE(5))\ne0$. We must observe that in this case $\gamma -2 = 8$.

\subsection{A non-stable bundle with $c_1 = c_2 = 0$, $\alpha = -4$, $\gamma = 9$.}
Take in $\PP^3$ homogeneous coordinates $x,y,z,t$ and let $Y$ be the non-reduced structure on the line $L\colon x=y=0$ defined by the ideal $I=(x^{8},x^{6}y^{2},x^{4}y^{4},x^{2}y^{6},y^{8},z^{6}x^{2}-y^{2}t^{6})$ (see \cite{CV1}, example~3.1, (iii)). Then $Y$ is $(-12)$-subcanonical and the zero locus of a section of a bundle $\shF = \shE(-4)$, where $\shE$ has $c_1 = 0$, $\alpha = -4$, $\gamma = 12$, $c_2 = 0$, $\delta = 16$.  By Theorem~\ref{non-stable}, iv), we see that $h^1(\shE(13))\ne0$, while the fact that $\gamma-2 = 10$ implies only that $h^1(\shE(10))\ne0$.

\subsection{Strictly semistable bundles with $c_1 = \alpha = 0$, $c_2 = 4$.}
Example $12$, $a$ of \cite{Ellia} has the following cohomology table:
\begin{center}
\renewcommand{\arraystretch}{1.25}
\begin{tabular}{|c|r r r r r r r|}
\hline
$h^1$ & 0 & 2 & 4 & 7 & 8 & 6 & 0 \\
$h^0$ & 0 & 0 & 0 & 1 & 4 & 10 & 20  \\
\hline
&  $-3$ & $-2$ & $-1$ & 0 & 1 & 2 & 3 \\
\hline
\end{tabular}
\end{center}
Observe that $\tau=\sqrt{6 c_2+1}-2 = \sqrt{24+1}-2 = 3$, hence Theorem~\ref{non-stable} gives a sharp bound when $c_1 = \alpha = 0$.

\begin{remark}
The above examples from 4.5 to 4.9 have also been checked with Macaulay 2 (see \cite{Macaulay}) which, of course, gives all the non-vanishing results that we obtain theoretically. It also shows that sometimes our results are not sharp. For instance, according to Macaulay 2, we have: $h^1(\shE(34))\ne0$ and $h^1(\shE(35)) = 0$ in example 4.7.
\end{remark}

\bigskip
\noindent
VALABREGA Paolo, Dipartimento di Matematica, Politecnico di Torino, 
Corso Duca degli Abruzzi 24, 10129 Torino, Italy,
e--mail: \texttt{paolo.valabrega@polito.it}

\noindent
VALENZANO Mario, Dipartimento di Matematica, Universit\`a di Torino, 
Via Carlo Alberto 10, 10123 Torino, Italy, 
e--mail: \texttt{mario.valenzano@unito.it}

\end{document}